\documentclass[12pt, reqno]{amsart}
\usepackage[T1]{fontenc}
\usepackage{dsfont}
\usepackage{mathrsfs}
\usepackage[colorlinks, citecolor=blue, linkcolor=blue]{hyperref}
\usepackage{xcolor}
\usepackage[a4paper,asymmetric]{geometry}
\usepackage{mathscinet}
\usepackage{latexsym}
\usepackage{amsthm}
\usepackage{amssymb}
\usepackage{amsfonts}
\usepackage{amsmath}
\usepackage{longtable}
\usepackage{graphicx}
\usepackage{multirow}
\usepackage{multicol}
\usepackage{amsfonts, amsmath}
\usepackage{latexsym,bm,amsfonts,amssymb,pifont,mathbbol,bbm}
\usepackage{verbatim}

\newtheorem{theorem}{Theorem}[section]
\newtheorem{thm}[theorem]{Theorem}
\newtheorem{lemma}[theorem]{Lemma}

\newtheorem{remark}[theorem]{Remark}
\newtheorem{proposition}[theorem]{Proposition}
\newtheorem{prop}[theorem]{Proposition}

\newtheorem{hyp}[theorem]{HYPOTHESIS}
\theoremstyle{definition}

\newtheorem{ex}[theorem]{Example}
 \newtheorem{nott}{Notation}
\theoremstyle{remark}
\numberwithin{equation}{section}

 \DeclareMathAlphabet{\mathpzc}{OT1}{pzc}{m}{it}
 \DeclareMathAlphabet{\mathsfsl}{OT1}{cmss}{m}{sl}

  \newcommand{\FH}{\mathfrak{H}}

\newcommand{\dif}{\mathrm{d}}

\newcommand{\R}{\mathbb{R}}

\newcommand{\N}{\mathbb{N}}
\newcommand{\Z}{\mathbb{Z}}

\newcommand{\abs}[1]{\left\vert#1\right\vert}
\newcommand{\set}[1]{\left\{#1\right\}}

\newcommand{\norm}[1]{\left\Vert#1\right\Vert}
 \newcommand{\innp}[1]{\langle {#1}\rangle}

\newcommand{\E}{\mathbb{E}}

\allowdisplaybreaks

\begin{document}
\title[B-E bounds and ASCLT for quadratic variation of Gaussian processes]{Berry-Ess\'{e}en bounds and almost sure CLT for the quadratic variation of a general Gaussian process}
\author[Y. Chen]{Yong CHEN}
 \address{School of Mathematics and Statistics, Jiangxi Normal University, Nanchang, 330022, Jiangxi, China}
\email{zhishi@pku.org.cn;\, chenyong77@gmail.com}
 \author[Z. Ding]{Zhen DING}
 \address{School of Mathematics and Statistics, Jiangxi Normal University, Nanchang, 330022, Jiangxi, China}
  \author[Y. Li]{Ying LI}
 \address{School of Mathematics and Computional Science, Xiangtan University, Xiangtan, 411105, Hunan, China. (Corresponding author.)}  \email{liying@xtu.edu.cn}
\begin{abstract}
In this paper, we consider the explicit bound for the second-order approximation of the quadratic variation of a general fractional Gaussian process $(G_t)_{t\ge 0}$. The second order mixed partial derivative of the covariance function $ R(t,\, s)=\mathbb{E}[G_t G_s]$ can be decomposed into two parts, one of which coincides with that of fractional Brownian motion and the other of which is bounded by $(ts)^{H-1}$ up to a constant factor.
This condition is valid for a class of continuous Gaussian processes that fails to be self-similar or have stationary increments. 
 Under this assumption, we obtain the optimal Berry-Ess\'{e}en bounds when $H\in (0,\,\frac23]$ and the upper Berry-Ess\'{e}en bounds when $H\in (\frac23,\,\frac34]$. As a by-product, we also show the almost sure central limit theorem (ASCLT) for  the quadratic variation when $H\in (0,\,\frac34]$. The results extend that of \cite{NP 09} to the case of general Gaussian processes, improve the Berry-Ess\'{e}en bounds  and unify the proofs in \cite{Tu 11}, \cite{AE 12} and \cite{KL 21} for respectively the sub-fractional Brownian motion, the bi-fractional Brownian motion and the sub-bifractional Brownian motion. \\



\end{abstract}
\maketitle

\section{ Introduction}\label{sec 03}
We are interested in the explicit bound for the second-order approximation of the quadratic variation of a general fractional Gaussian process $G=\set{G_t:\,t\ge 0}$ on $[0,T]$, defined as
\begin{align*} 
Z_n&:=\sum_{k=0}^{n-1}\Big[  (G_{k+1}-G_{k })^2- \E[(G_{k+1}-G_{k})^2 ] \Big].
\end{align*}
Let the renormalization of $Z_n$ be 
\begin{align} \label{vn defn}
V_n&:= \frac{1}{\sigma_n}Z_n,
\end{align}
where $\sigma_n>0$ is so that $\E[V_n^2]=1$, i.e., $\sigma_n^2:= \E[Z_n^2]$. 

When $G$ is the fractional Brownian motion $B^{H}$, by using Stein method and
Malliavin calculus, Nourdin and Peccati \cite{ND 12, NP 09} derived explicit bounds for the total variation distance between the law of  $V_n$ and the standard normal law $N$. From then on, the same problem was extended to some other fractional Gaussian processes such as the sub-fractional Brownian motion, the bi-fractional Brownian motion and the sub-bifractional Brownian motion in \cite{Tu 11}, \cite{AE 12} and \cite{KL 21} respectively. 

We find that the above four types of fractional Gaussian processes are all special examples of the following general Gaussian process $G$. Let \begin{equation} \label{exmp1}
 R^{B}(s,t)= \frac{1}{2}(\abs{s}^{2H} + |t|^{2H} - |t-s|^{2H}),
 \end{equation} be the covariance function of the fractional Brownian motion $\{B^H(t), t\geq 0\}$. 
\begin{hyp}\label{hypthe 1} 
 For $H \in (0,\,1)$,  the covariance function $R(t,s)=\mathbb{E}[G_{t}G_{s}]$ satisfies that
 \begin{enumerate}
	\item  for any $s\geq 0$, $R(0,s)=0$.
	\item  for any fixed $s\in (0,T)$, $ R(t,s) $ is continuous on $[0,T]$ and differentiable function with respect to $t$ on $(0,s)\cup(s,T)$, such that $\frac{\partial  }{\partial t }R(t,s) $ is absolutely integrable.
	\item  
for any fixed $t\in (0,T)$, the difference 
$$\frac{\partial R(t,s)}{\partial t} - \frac{\partial R^{B}(t,s)}{\partial t} $$
is continuous on $[0,T]$, and it is differentiable with respect to $s$ on $(0,T)$ such that $\Psi (t,s)$, the partial derivative with respect to $s$ of the difference, satisfies \begin{align}\label{cond hyp2}
    |\Psi (t,s)|&\leq C_{H }^{^{\prime }}|ts|^{H -1},
\end{align}
where the constants $H$,
 $C_{H }^{^{\prime }}\geq 0$ do not depend on T, and $R^B(t,s)$ is the covariance function of the fractional Brownian motion as in \eqref{exmp1}.
	 \end{enumerate}
\end{hyp}
\begin{ex}
The subfractional Brownian motion $\{S^H(t), t \geq 0\}$ with parameter $H\in (0,1)$ has the covariance function
$$R(t,s)=s^{2H}+t^{2H}-\frac{1}{2}\left((s+t)^{2H}+|t-s|^{2H}\right),$$
which satisfies Hypothesis~\ref{hypthe 1}.
\end{ex}
\begin{ex}
The bi-fractional Brownian motion $\{B^{H',K}(t), t\geq 0\}$ with parameters $H',K \in (0, 1)$ has the covariance function
$$R(t,s)=\frac{1}{2}\left((s^{2H'}+t^{2H'})^K - |t-s|^{2H'K}\right) ,$$
which satisfies Hypothesis~\ref{hypthe 1} when $H := H'K $. 
\end{ex}

\begin{ex}\label{exmp5}
The generalized sub-fractional Brownian motion $S^{H',K}(t) $ with parameters $H' \in (0, 1),\,K \in[1,2)$ and $H'K\in (0,1)$ satisfies Hypothesis~\ref{hypthe 1} when $H := H'K $. The covariance function 
is 
$$ R(t,\, s)= (s^{2H'}+t^{2H'})^{K}-\frac12 \big[(t+s)^{2H'K} + \abs{t-s}^{2H'K} \big].$$
\end{ex}


 
 \begin{nott}
 Given two deterministic numeric sequences $(a_n)_{n\geq 0}$,$(b_n)_{n\geq 0}$, we use the following notations and definitions for respectively commensurability, equivalence:
\begin{align*}
a_n&\asymp b_n\iff \exists c,C>0:cb_n\geq a_n\leq Cb_n,\,\text{for $n$ large enough,}\\
a_n&\sim b_n \iff \exists c_n,C_n>0:\lim_{n\to \infty}c_n=\lim_{n\to \infty}C_n=1, \,  \, c_nb_n\leq a_n\leq C_nb_n,\,\text{for $n$ large enough}.
\end{align*}
\end{nott}
 Now we list our main results as follows:
 \begin{thm}\label{main thm 1}
Let $N\sim{N(0,1)}$ and $V_n$ be given as in \eqref{vn defn} and suppose Hypothesis~\ref{hypthe 1} holds.  When  $H\in (0,\,\frac23)$, 
 \begin{equation}\label{d}
  d_{TV}(V_n,\,N) \asymp  n^{-\frac12}; \end{equation}
when $H=\frac23$, 
 \begin{equation}\label{d 1}
  d_{TV}(V_n,\,N) \asymp  n^{-\frac12}\log^2 n ; \end{equation}
when $H\in ( \frac23,\,\frac{3}{4})$, there exists  a positive contant $c_{H}$ depending on $H$ such that for any $n\geq1$,
 \begin{equation}\label{d 2}
d_{TV}(V_n,\,N)\leq c_{H }\times  n^{\frac12 (4H-3)} ;
 \end{equation}
 when $H=\frac{3}{4}$, there exsits  a positive contant $c $  such that for any $n\geq1$,
  \begin{equation}\label{d 3}
d_{TV}( {V_n}, \,N)\leq   \frac{c}{ (\log n)^{\frac32}} ;
 \end{equation}
 \end{thm}
 \begin{remark}  \begin{enumerate}
	\item  The above Berry-Ess\'{e}en types bounds are more sharp than those obtained in \cite{Tu 11}, \cite{AE 12} and \cite{KL 21} for respectively the sub-fractional Brownian motion, the bi-fractional Brownian motion and the sub-bifractional Brownian motion.
	\item We do not know how to obtain the optimal bound in the case of $H\in ( \frac23,\,\frac{3}{4}]$.	
	\item In the same vein, we can also extend Theorem~\ref{main thm 1} to the $p$th Hermite variation with $p>2$. 
\end{enumerate}
	\end{remark} 

As a by-product of Theorem~\ref{main thm 1}, we have the ASCLT of the sequence $(V_n)_{n\geq 1}$.
 \begin{thm}\label{main thm 2}
If $H\in (0,\frac34]$ then the sequence $(V_n)_{n\geq 1}$ satisfies the ASCLT. In other words, for any bounded and continuous function $\varphi: \R {\rightarrow} \R$, we have almost surely,
\begin{equation}
\frac{1}{\log n}\sum_{k=1}^{n}\frac1k\varphi(V_k)\stackrel{a.s.}{\longrightarrow}\E\varphi(N)
\end{equation}
as $n{\rightarrow}\infty$, where $N\sim{N(0,1)}$.
 \end{thm}

\section{Preliminary}
In this section, we describe some basic facts on stochastic calculus with respect to the Gaussian process, 
for more complete presentation on the subject can be found in \cite{ Jolis2007}.\par
Denote $G=\set{G_{t},t\in[0,T]}$ as a continuous centered Gaussian process with $G_0=0$ and the covariance function 
\begin{equation}
\mathbb{E}(G_{t}G_{s})=R(s,t),\,\, s,t\in[0,T],
\end{equation}defined on a complete probability space$(\Omega,\mathcal{F},P)$, where the $\mathcal{F}$ is generated by the Gaussian family $G$. 
Suppose in addition that the covariance function $R$ is continuous.
Let $\mathcal{E}$ denote the space of all real valued step functions on $[0,T]$. The Hilbert space $\mathfrak{H}$ is defined as the closure of $\mathcal{E}$ endowed with the inner product
\begin{equation}\label{innp qishi}
\langle\mathbbm{1}_{[a,b)},\mathbbm{1}_{[c,d)}\rangle_{\mathfrak{H}}=\mathbb{E}((G_{b}-G_{a})(G_{d}-G_{c})),
\end{equation} where $\mathbbm{1}_{[a,b)}$ is the indicator function of the interval $[a,b)$.
With abuse of notation, we also denote $G=\set{G(h), h\in{\mathfrak{H}}}$ as the isonormal Gaussian process on the probability space, indexed by the elements in the Hilbert space $\mathfrak{H}$. Then $G$ is a Gaussian family of random variables such that
\begin{equation}
\mathbb{E}(G)=\mathbb{E}(G(h))=0,\,\,
\mathbb{E}(G(g)G(h))=\langle g,h\rangle_{\mathfrak{H}},\quad \forall g,h \in\mathfrak{H}.
\end{equation}
\begin{nott}\label{innp h1 h2}
Let $R^B(t,s)$ be the covariance function of the fractional Brownian motion as in \eqref{exmp1}. $\mathcal{V}_{[0,T]}$ denote  the set of bounded variation functions on $[0,T]$. For  functions $f,\,g \in\mathcal{V}_{[0,T]}$, we define two  products as
\begin{equation}
\begin{aligned}
\innp{f,\,g }_{\mathfrak{H}_1}=- \int_{[0,T]^{2}}f(t ) \frac{\partial R^B(t,s)}{\partial t}  \dif t \nu_{g}({\dif s}),\\
\innp{f,\,g }_{\mathfrak{H}_2}=C_{H}^{'} \int_{[0,T]^{2}} \abs{f(t) g(s)}(ts)^{H-1}\dif t  \dif s.
\end{aligned}
\end{equation}
\end{nott}
The following proposition is an extension of \cite[Theorem 2.3]{Jolis2007} and \cite[Proposition 2.2]{CL 21}, which gives the inner products representation of the Hilbert space $\mathfrak{H}$:
\begin{proposition}\label{inner products represent}
$\mathcal{V}_{[0,T]}$ is dense in $\mathfrak{H}$ and we have
  \begin{align}  \label{inner product 001}
\innp{f,g}_{\FH}=\int_{[0,T]^2} R(t,s) \nu_f( \dif t)   \nu_{g}( \dif s),\qquad \forall f,\, g\in \mathcal{V}_{[0,T]},  
\end{align}where $\nu_{g}$ is the restriction to $([0,T ], \mathcal{B}([0,T ]))$ of 
the Lebesgue-Stieljes signed measure associated with $g^0$ defined as
\begin{equation*}
g^0(x)=\left\{
      \begin{array}{ll}
 g(x), & \quad \text{if } x\in [0,T).\\
0, &\quad \text{otherwise}.     
 \end{array}
\right.
\end{equation*}
Furthermore, if the covariance function $R(t,s)$ satisfies Hypothesis~\ref{hypthe 1}, then  
  \begin{align} \label{innp fg3}
\innp{f,g}_{\FH}=-\int_{[0,T]^2}  f(t) \frac{\partial R(t,s)}{\partial t} \dif t  \nu_{g}(\dif s),\qquad \forall f,\, g\in \mathcal{V}_{[0,T]}.
\end{align}
In particular, 
 we have
\begin{equation}\label{inequality 29}
\abs{\innp{f,\,g}_{\mathfrak{H}} - \innp{f,\,g }_{\mathfrak{H}_1}}\leq  {\innp{f,\,g }_{\mathfrak{H}_2} },\qquad \forall f,\, g\in \mathcal{V}_{[0,T]}.
\end{equation}
\end{proposition}
\begin{remark}
 When $H\in (\frac12,\,1)$,  Hypothesis~\ref{hypthe 1} (3) and Lemma~\ref{lemm integra by part} imply that the identity \eqref{innp fg3} can be rewritten as
\begin{equation}\label{inner product 000}
\langle f,g\rangle_\mathfrak{H}=\int_{[0,T]^{2}}f(t)g(s)\frac{\partial ^{2}R(t,s)}{\partial
t\partial s}dtds,          \quad        \forall f,g\in \mathcal{V}_{[0,T]}.
\end{equation} In this case, the inequality \eqref{inequality 29} has obtained from \eqref{inner product 000} in \cite{CZ 21}. 
 When $H\in (0,\frac12)$, it is well known that both $\frac{\partial^2 }{\partial t \partial s}R(t,s)$ and $\frac{\partial^2 }{\partial t \partial s}R^B(t,s)$ are not absolutely integrable. But the absolute integrability of their difference makes the key inequality \eqref{inequality 29} still valid. 
\end{remark}
\begin{proof}
The first claim and the identity (\ref{inner product 001}) are rephrased from Theorem 2.3 of \cite{Jolis2007}. 
Hypothesis~\ref{hypthe 1} (2) and Lemma~\ref{lemm integra by part} imply the inner products representation \eqref{innp fg3}.

Finally, it follows from the identity \eqref{innp fg3} and Notation~\ref{innp h1 h2} that
\begin{align*}
\innp{f,\,g }_{\mathfrak{H}_1}-\innp{f,\,g}_{\mathfrak{H}}&=\int_{[0,T]^2}  f(t) 
\Big[ \frac{\partial R(t,s)}{\partial t} - \frac{\partial R^{B}(t,s)}{\partial t}\Big]\dif t  \nu_{g}(\dif s).
\end{align*}
By the fundamental theorem of calculus (see Proposition 1.6.41 of \cite{Tao 11}), Hypothesis~\ref{hypthe 1} (1) and (3)  imply that when $s\neq t$,
\begin{align}\label{difference r rb}
\frac{\partial R(t,s)}{\partial t} - \frac{\partial R^{B}(t,s)}{\partial t}=\int_{0}^s \Psi(t,r)\dif r.
\end{align}
Hence, Lemma~\ref{lemm integra by part} implies that
\begin{align*}
\innp{f,\,g}_{\mathfrak{H}}-\innp{f,\,g }_{\mathfrak{H}_1}&=-\int_{[0,T] }  f(t) \dif t\int_{[0,T]}
\Big[ \frac{\partial R(t,s)}{\partial t} - \frac{\partial R^{B}(t,s)}{\partial t}\Big]  \nu_{g}(\dif s)\\
&=\int_{[0,T] }  f(t) \dif t\int_{[0,T]} g(s)
\Psi(t,s)  \dif s,
\end{align*} which implies the inequality \eqref{inequality 29} since $\Psi(t,s) $  satisfies the inequality \eqref{cond hyp2}.
\end{proof}

Denote $\mathfrak{H}^{\otimes p}$ and $\mathfrak{H}^{\odot p}$ as the $p$th tensor product and the $p$th symmetric tensor product of the Hilbert space $\mathfrak{H}$. Let $\mathcal{H}_p$ be the $p$th Wiener chaos with respect to $G$. It is defined as the closed linear subspace of $L^2(\Omega)$ generated by the random variables $\{H_p(G(h)): h \in \mathfrak{H}, \ \|h\|_{\mathfrak{H}} = 1\}$, where $H_p$ is the $p$th Hermite polynomial defined by
$$H_p(x)= {(-1)^p}  {\rm e}^{\frac{x^2}{2}} \frac{d^p}{dx^p} \rm{e}^{-\frac{x^2}{2}}, \quad p \geq 1,$$
and $H_0(x)=1$. We have the identity $I_p(h^{\otimes p})=H_p(G(h))$ for any $h \in \mathfrak{H}$ with $\norm{h}_{\mathfrak{H}} = 1 $ where $I_p(\cdot)$ is the $p$th multiple Wiener-It\^o integral. Then the map $I_p$ provides a linear isometry between $\mathfrak{H}^{\odot p}$ (equipped with the norm $ {\sqrt{p!}}\|\cdot\|_{\mathfrak{H}^{\otimes p}}$) and $\mathcal{H}_p$. Here $\mathcal{H}_0 = \mathbb{R}$ and $I_0(x)=x$ by convention.

The following Theorem \ref{fm.theorem}, known as the optimal fourth moment theorem, provides exact rates of convergence in total variation distance between a multiple Wiener-It\^{o} integral and a normal distribution (see \cite{NP 15, BBNP 12}).
\begin{thm} \label{fm.theorem}
 Let $N \sim N(0, 1)$ be a standard Gaussian random variable. Let $\set{F_n:\,n\ge 1}$ be a sequence of random variables living in the $p$th multiple Wiener-It\^o integral with unit variance. If  $\lim_{n \to \infty} \mathbb{E}[F_n^4] = { 3}$, then there exist two finite constants $0 < c < C$ (possibly depending on $p$ and on the sequence $\set{F_n}$, but not on $n$) such that the following estimate in total variation holds for
every $n$:
 \begin{equation*}
 c\mathbf{M}(F_n)\le d_{TV}(F_n,\, N)\le C \mathbf{M}(F_n),
 \end{equation*}where 
 \begin{equation*}
 \mathbf{M}(F_n):=\max\set{\abs{\E[F_n^3]}, \E[F_n^4-3]}.
 \end{equation*}
\end{thm}
The quantities $\kappa_3(F_n):= \E({F_n}^3)$ and $\kappa_4(F_n):= \E[{F_n}^4]-3$ are called the 3rd and 4th cumulants of $F_n$. That $\kappa_3(F_n)$ coincides with the third moment is because $F_n$ is centered. Moreover, $\kappa_4(F_n)$ is strictly positive (see \cite{NP 12, NV 16} ).

The following theorem is used to show the ASCLT.
\begin{theorem}[\cite{BNT 10}] \label{asclt chufad}
Let $p \geq2$ be an integer, and let $F_n=I_q(f_n)$, with $f_n\in{\FH}^{\odot p}$. Assume that for all n, and that $F_n\stackrel{ {law}}{\to}\mathcal{N}(0,1)$ as $n \to \infty$. If the two following conditions are satisfied
\begin{equation}\begin{aligned}
(1).\,\,  &\sum_{n=2}^{\infty}\frac1{n\log^{2} n}\sum_{k=1}^{n}\frac{1}{k}\norm{f_k\otimes_{r} f_k}_{{\FH}^{\otimes 2(q-r)}}<\infty,\,\, \text{for every}\;1\leq r\leq p-1,\\
(2).\,\, &\sum_{n=2}^{\infty}\frac1{n\log^{3} n}\sum_{k,l=1}^{n}\frac{|\langle f_k,f_l\rangle_{{\FH}^{\otimes p}}|}{kl}<\infty.
\end{aligned}\end{equation}
Then $ \set{F_n:\,n\ge 1}$ satisfies an ASCLT.
\end{theorem}

Denote
\begin{align*}
 \rho(r)&:=\frac12 \big( \abs{r+1}^{2H} +  \abs{r-1}^{2H}-2 \abs{r}^{2H}  \big), \quad r\in \Z. 
 \end{align*}
It is well-known that for any $r\neq 0$, $ \rho(r)$ is positive when $H\in (\frac12,\,1)$ and is negative when $H\in (0,\,\frac12)$. 
It behaves asymptotically as $\abs{r}\to \infty$, \begin{align*}
 \rho(r)\sim H(2H-1) \abs{r}^{2H-2}.
\end{align*} In particular, when $H>\frac12$, for $|r|$ large enough,
\begin{equation*}
\rho (r)\geq H(H-\frac12)(1+|r|)^{2H-2}.
\end{equation*}

The following proposition is cited from \cite[p.74]{ND 12}.
\begin{prop}\label{prop rho2 lim}
If $H\in (0,\,\frac34)$ then
\begin{equation}\label{sigma2 defn}
\lim_{n\to \infty} \frac{2}{n}\sum_{i,j=0}^{n-1} \rho^2(i-j)=2\sum_{r\in \Z} \rho^2(r):=\sigma^2;
\end{equation} and if $H=\frac34$ then 
\begin{equation*}
\lim_{n\to \infty} \frac{2}{n \log n}\sum_{i,j=0}^{n-1} \rho^2(i-j)=\frac{9}{16}.
\end{equation*} 
\end{prop}
The following propositions are cited from Proposition 6.6 and Proposition 6.7 of \cite{BBNP 12} respectively. For the case of $H=\frac34$, please refer to \cite{NV 16}.
\begin{prop}\label{3 jichu bijiao}
We have
 \begin{equation*}
\sum_{j,k,l=0}^{n-1}\rho(j-k)\rho(k-l)\rho(j-l)\asymp \left\{
      \begin{array}{ll}
n, & \quad \text{if }  H\in (0\,\frac23),\\
n\log^{2}n, & \quad \text{if } H = \frac23, \\
n^{6H-3}, &\quad \text{if } H\in (\frac23,\, \frac34].
  \end{array}
\right.
\end{equation*}
and
 \begin{equation*}
\sum_{i,j,k,l=0}^{n-1}\rho(i-j)\rho(i-k)\rho(k-l)\rho(j-l)\asymp \left\{
      \begin{array}{ll}
 n, & \quad \text{if }  H\in (0,\,\frac58),\\
 {n(\log n)^3,} & \quad \text{if } H = \frac58, \\
n^{8H-4}, &\quad \text{if } H\in (\frac58,\, \frac34].
  \end{array}
\right.
\end{equation*}
\end{prop}

\section{proof of Theorems~\ref{main thm 1} and \ref{main thm 2}.}
 We will discuss exclusively the case $H \neq \frac12$ since the case $H= \frac12$ is easy. First, we need two technical lemmas. The first one is a variant of the classical integration by parts formula. 
  It is well known that when $f,\phi:\,\R\to \R$ are monotone non-decreasing and continuous functions, then 
\begin{equation}\label{tao formula}
-\int_{[a,b]}  f\dif \phi= \int_{[a,b]}  \phi \dif f+f(a)\phi(a) - f(b)\phi(b)
\end{equation} for any compact interval $[a,b]$  \cite[p.160]{Tao 11}. 
If $f $ are continuously differentiable on $[a,b]$, we formally have the following expression (see (10) of \cite{CL 21}):  
\begin{align*}
\dif \big(f\cdot \mathbf{1}_{[a,b)}(\cdot)\big)=\Big[f'(t) \mathbf{1}_{[a,b]}(t) +f(t)\big(\delta_a(t)-\delta_b(t)\big)\Big] \dif t,
\end{align*}
This expression suggests that we can take the terms $f(a)\phi(a) - f(b)\phi(b)$ in the right hand side of \eqref{tao formula} as parts of the measure $\nu_{f}$ for conveniences. This is what  the following lemma to do.
 
\begin{lemma}(Integration by parts formula)\label{lemm integra by part}
 Let $[a,b] $ be a compact interval of positive length, let $\phi:[a,b]\to \R$ be continuous on $[a,b]$ and differentiable in $(a,b)$, such that $\phi'  $ is absolutely integrable. For any  $f\in \mathcal{V}_{[a,b]}$,  we have
\begin{align}\label{integ by parts 001}
-\int_{[a,b]}  f(t) \phi'(t) \dif t =\int_{[a,b]} \phi(t) \nu_f( \dif t),
\end{align}where $\nu_f$ is given as in Proposition~\ref{inner products represent}, i.e., $\nu_{f}$ is the restriction to $([a,b ], \mathcal{B}([a,b ]))$ of 
the Lebesgue-Stieljes signed measure associated with $f^0$ defined as
\begin{equation*}
f^0(x)=\left\{
      \begin{array}{ll}
f(x), & \quad \text{if } x\in [a,b).\\
0, &\quad \text{otherwise }.     
 \end{array}
\right.
\end{equation*}
\end{lemma}
\begin{proof} The proof is similar to that of \cite[Proposition 2.2]{CL 21}. We establish this in stages.
We first deal with the case when $f$ is a step functions on $[a,b )$ of the form
\begin{align*}
f=\sum_{j=0}^{N-1}f_j \mathbf{1}_{[t_j,t_{j+1})},
\end{align*}
where $\set{a= t_0 < t_1 < \cdots < t_N = b } $ is a partition of $[a,b )$ and $f_j\in\R$. The corresponding signed measure is \cite[p.1123]{Jolis2007}
\begin{align*}
\nu_{f}=\sum_{j=1}^{N-1}( f_j-f_{j-1})\delta_{t_j}+f(a+)\delta_0-f(b-)\delta_T.
\end{align*} By the fundamental theorem of calculus again, the assumption of $\phi$ implies that 
\begin{align*}
\int_{[t_j,\,t_{j+1})}  \phi'(t)  \dif t = \phi(t_{j+1})-\phi (t_j). 
\end{align*}
Hence,  the following formula of integration by parts hold:  
\begin{align} 
-\int_{[a,b]}  f(t)  \phi'(t) \dif t &=\sum_{j=0}^{N-1} f_j \big[\phi (t_j)- \phi(t_{j+1})\big]\nonumber \\
&=\sum_{j=1}^{N-1}( f_j-f_{j-1})\phi({t_j})+f_0\phi(t_0)-f_{N-1}\phi(t_{N})\nonumber \\
&=\int_{[a,b]} \phi(t) \nu_f( \dif t).\label{integ by parts 1}
\end{align}

Now we assume that $f$ is a right continuous monotone non-decreasing function on $[a, b)$. For any sequence partitions $\pi_n= \set{a = t_0^n < t_1^n < \cdots < t_{k_n}^n = b } $ such that $\pi_n\subset \pi_{n+1}$ and $\abs{\pi_n}\to 0$ as $n \to \infty$, consider 
\begin{align*}
f_n=\sum_{j=0}^{{k_n}-1} f(t_j^n) \mathbf{1}_{[t_j^n,t_{j+1}^n)},
\end{align*}which is uniform bounded since $f$ is bounded.
It is clear that the sequence of signed measures $\nu_{f_n}$ converges weakly to $\nu_{f}$ \cite{Jolis2007}. Hence,  we have
\begin{align*}
\int_{[a,b]} \phi(t) \nu_f( \dif t)&=\lim_{n\to \infty} \int_{[a,b]} \phi(t) \nu_{f_n}( \dif t) \\
&=-\lim_{n\to \infty} \int_{[a,b]}  f_n(t)  \phi'(t) \dif t\qquad \text{(by \eqref{integ by parts 1})}\\
&=-\int_{[a,b]}  f(t)  \phi'(t) \dif t,
\end{align*}where the last line is from Lebesgue's dominated theorem since $f_n$ is uniform bounded and $\phi'$ are absolutely integrable.

Finally, it is well known that every function of bounded variation is the difference of two bounded monotone non-decreasing function and that the value of $f$ at its points of discontinuity are irrelevant for the purposes of determining  the Lebesgue-Stieltjes measure $\nu_f$ \cite{Tao 11}. Hence, (\ref{integ by parts 001}) is valid for any  $f\in \mathcal{V}_{[a,b]}$.
\end{proof}

\begin{lemma} \label{beta function inequality}
Let $ v_1,\cdots, v_l$ be positive.
There is a positive constant $c$ depending on $v_1,\cdots, v_l$ such that when $r\in \N:=\set{1,2,\cdots}$ is large enough,
\begin{align}\label{dirichlet inequ 1}
\sum_{r_i\in \N, \sum_{i=1}^{l} r_i< r } r_1^{v_1-1}r_2^{v_2-1}\cdots r_l^{v_l-1} \le c \times r^{ \sum_{i=1}^l v_i}.
\end{align}
\end{lemma}
\begin{remark}
When $ v_1,\cdots, v_l$ are negative, the following inequality is trivial: there is a positive constant $c$ depending on $v_1,\cdots, v_l$ such that when $r\in \N $ is large enough,
\begin{align}\label{dirichlet inequ 001}
\sum_{r_i\in \N, \sum_{i=1}^{l} r_i< r } r_1^{v_1-1}r_2^{v_2-1}\cdots r_l^{v_l-1} \le c <\infty.
\end{align}

\end{remark}
\begin{proof} Let $v_0>0$.
It is well-known that on  the standard simplex in $\R^{l}$: $$ {T}^{l}:=\set{(x_1,\cdots, x_l): \,x_i \ge 0,\sum_{i=1}^{l} x_i \le1},$$  the following integral converges:
\begin{align*}
\int_{T^l} x_1^{v_1-1}x_2^{v_2-1}\cdots x_l^{v_l-1}(1-x_1-\cdots-x_l)^{v_0-1}\dif x=\frac{\Pi_{i=0}^l\Gamma(v_i)}{\Gamma(\sum_{i=0}^l v_i)}
\end{align*} where $\Gamma(\cdot)$ the Gamma function.
The change of variables implies that there is a positive constant $c$ depending on $v_0,v_1,\cdots, v_l$ such that when $r\in \N$ is large enough,
\begin{align*}
\sum_{r_i\in \N, \sum_{i=1}^{l} r_i< r } r_1^{v_1-1}r_2^{v_2-1}\cdots r_l^{v_l-1}(r-r_1-\cdots-r_l)^{v_0-1}\le c \times r^{-1+\sum_{i=0}^l v_i}.
\end{align*}
Especially, when $v_0=1$, the above inequality collapses to \eqref{dirichlet inequ 1}.
\end{proof}
Without any loss of generality, we suppose for simplicity that 
 $C_{H}'=1 $ in this section.
Denote
\begin{align}
\theta({i,j})&:=\E\Big[\big(G_{i+1}-G_i\big)\big(G_{j+1}-G_j\big)\Big], \nonumber\\
\gamma({i,j})&:= \theta(i,j)-\rho(i-j).\label{gammaij dfn}
\end{align}
\begin{prop}\label{prop gamma2 lim}
Under Hypothesis~\ref{hypthe 1},  there exists  a positive constant $c$ such that for any $n\ge 1$,
\begin{align}\label{gamma2 div n estimate}
\frac{1}{n}\sum_{i,j=0}^{n-1} \gamma(i,j)^2\le c\times n^{(4H -3)\vee (-1)} .
\end{align}Hence, we have that as $n\to \infty$,  when $H\in (0,\,\frac34)$,
\begin{align*}
\frac{1}{n}\sum_{i,j=0}^{n-1} \gamma(i,j)^2 \to 0;
\end{align*}and  when $H=\frac34$,
\begin{align*}
\frac{1}{n\log n}\sum_{i,j=0}^{n-1} \gamma(i,j)^2 \to 0.
\end{align*}
\end{prop}
\begin{proof} It is clear that we need only show the inequality \eqref{gamma2 div n estimate} holds. It follows \eqref{innp qishi}, the definition of the inner product, that 
\begin{align*}
\theta(i,j)&=\E\Big[\big(G_{i+1}-G_i\big)\big(G_{j+1}-G_j\big)\Big]=\innp{\mathbbm{1}_{[i,i+1)},\,\mathbbm{1}_{[j,j+1)}}_{\FH},
\end{align*} and 
\begin{align*}
\rho(i-j)&=\E\Big[\big(B^{H}_{i+1}-B^{H}_i\big)\big(B^{H}_{j+1}-B^{H}_j\big)\Big]=\innp{\mathbbm{1}_{[i,i+1)},\,\mathbbm{1}_{[j,j+1)}}_{\FH_1},
\end{align*} where $B^{H}$ is the fractional Brownian motion with Hurst index $H$. The inequality \eqref{inequality 29} implies that 
\begin{align}
\abs{ \gamma(i,j)}=\abs{ \theta(i,j)-\rho(i-j)}&\le \int_{[0,T]^{2}} { \mathbbm{1}_{[i,i+1)}(r_{1}) \mathbbm{1}_{[j,j+1)}(r_{2})}(r_{1}r_{2})^{H-1}dr_{1}dr_{2}\nonumber \\
&=\frac{1}{H^2}\big[(i+1)^{H}- i^{H}\big]\big[(j+1)^{H}- j^{H}\big].\label{gamma estimate}
\end{align}Hence,  as $n\to \infty$, 
\begin{align}
\frac{1}{n}\sum_{i,j=0}^{n-1} \gamma(i,j)^2&\leq \frac{1}{n H^4}\sum_{i,j=0}^{n-1}  \big[(i+1)^{H}- i^{H}\big]^2\big[(j+1)^{H}- j^{H}\big]^2\nonumber\\
&=\frac{1}{n } \Big(\frac{1}{H^2} \sum_{i=0}^{n-1}\big[(i+1)^{H}- i^{H}\big]^2\Big)^2 .\label{gammaij estimate 1}
\end{align}
It is clear that the function $f(u)=(1+u)^{H}$ in $u\in [0,\infty)$ is concave, i.e., $f''(u)\le 0$, which implies that for any $u\ge 0$, $(1+u)^{H}\le 1+H u$. Hence, for any $i\ge 1$,
\begin{equation}\label{gamma beta inquality}
(i+1)^{H}- i^{H}= i^{H} \big[(1+\frac{1}{i})^{H}-1\big]\le H \times  i^{H-1}.
\end{equation}
 We have that there exists a positive constant $c$ independent on $ n$ such that
\begin{align}\label{gamma estimate 2}
\frac{1}{H^2} \sum_{i=0}^{n-1}\big[(i+1)^{H}- i^{H}\big]^2\le \frac{1}{H^2}+\sum_{i=1}^{n-1} i^{2(H-1)}\le c\times  {n^{(2H-1)\vee 0}},
\end{align} please refer to Lemma 6.3 of \cite{ND 12}.

By plugging the above inequality into \eqref{gammaij estimate 1}, we obtain the desired inequality \eqref{gamma2 div n estimate}. 
\end{proof}

\begin{prop} \label{sigma2 limit prop}
Recall that $\sigma_n^2:=\E[Z_n^2]$ and $\sigma^2$ is given as in \eqref{sigma2 defn}.
Under Hypothesis~\ref{hypthe 1},  we have  \\
(i) when $H\in (0,\,\frac34)$,  as $n\to \infty$,
\begin{align}\label{sigman 2 limit1}
\frac{\sigma_n^2}{n} \to \sigma^2.
\end{align}
(ii) when $H=\frac34$,  as $n\to \infty$,
\begin{align}\label{sigman 2 limit12}
\frac{\sigma_n^2}{n \log n} \to \frac{9}{16}.
\end{align}
\end{prop}
\begin{proof}
By the definition of multiple Wiener-It\^o integrals, we rewrite $Z_n$ as follows:
 \begin{equation}\label{zn exp}
Z_{n}  =I_{2}(g_{n}),
\end{equation} 
where \begin{equation}\label{gn exp}
g_{n}= \sum_{i=0}^{n-1} \mathbbm{1}_{[i,i+1)}^{\otimes 2}.\end{equation}

By It\^o's isometry, we have 
\begin{align}
\sigma^2_n= \E[Z_n^2]&=2 \norm{g_n}_{\FH}^2=2\sum_{i,j=0}^{n-1} \innp{\mathbbm{1}_{[i,i+1)},\,\mathbbm{1}_{[j,j+1)} }_{\FH}^2= 2\sum_{i,j=0}^{n-1} \theta(i,j)^2.
\end{align}
It is clear that the identity \eqref{gammaij dfn} implies that
\begin{align}\label{prop theta2}
\frac{1}{n}\sum_{i,j=0}^{n-1} \theta(i,j)^2&=  \frac{1}{n}\sum_{i,j=0}^{n-1} \rho^2(i-j) + \frac{1}{n}\sum_{i,j=0}^{n-1} \gamma^2(i, j) + \frac{2}{n}\sum_{i,j=0}^{n-1} \rho(i-j)\gamma(i, j) .\end{align}
The Cauchy-Schwarz inequality implies that the third term satisfies that as $n\to \infty$,
\begin{equation*}
\abs{\frac{2}{n}\sum_{i,j=0}^{n-1} \rho(i-j)\gamma({i, j})}\le 2 \Big(\frac{1}{n}\sum_{i,j=0}^{n-1} \rho^2(i-j) \times \frac{1}{n}\sum_{i,j=0}^{n-1} \gamma({i, j})^2 \Big)^{\frac12}\to 0,
\end{equation*}where in the last line we have used Propositions \ref{prop rho2 lim} and \ref{prop gamma2 lim}. 
By plugging this limit into the identity \eqref{prop theta2} and using Propositions \ref{prop rho2 lim} and \ref{prop gamma2 lim} again, we obtain the desired limit \eqref{sigman 2 limit1}. 

In the similar vein, the desired limit \eqref{sigman 2 limit12} holds.
\end{proof}

\begin{prop}\label{rop 3 bijiao}
Let $\theta(i,j),\,\gamma(i, j),\,\rho(r )$ be given as in \eqref{gammaij dfn}. When $H\in (0,\,1)$, there exists a positive constant $c$ such that for all $n\ge 1$
\begin{align}
\abs{\sum_{j, k, l=0}^{n-1}\gamma(j,k)\gamma(k,l)\gamma(j,l)}&\leq c\times n^{(6H-3)\vee 0}, \nonumber \\ 
\abs{\sum_{j, k, l=0}^{n-1}\gamma(j,k)\gamma(k,l)\rho(j-l)}&\leq c\times n^{(6H-3)\vee 0},\label{2 eiitmate main}\\
\abs{\sum_{j, k, l=0}^{n-1}\gamma({j,k})\rho(k-l)\rho(j-l)}& \leq c\times n^{(6H-3)\vee 0}.\label{3 eiitmate main}
\end{align}
\end{prop}
\begin{proof}It follows from the inequalites \eqref{gamma estimate} and \eqref{gamma estimate 2} that 
\begin{align*}
\abs{\sum_{j, k, l=0}^{n-1}\gamma(j,k)\gamma(k,l)\gamma(j,l)}&\leq \sum_{j, k, l=0}^{n-1}\abs{\gamma(j,k)\gamma(k,l)\gamma(j,l)}\\
&\leq \Big( \frac{1}{H^2} \sum_{j=0}^{n-1}\big[(j+1)^{H}- j^{H}\big]^2\Big)^3\\
&\le c\times n^{(6H-3)\vee 0}.
\end{align*}
Similarly, we have 
\begin{align}
&\abs{\sum_{j, k, l=0}^{n-1}\gamma(j,k)\gamma(k,l)\rho(j-l)}\nonumber \\
&\leq \sum_{j, k, l=0}^{n-1}\abs{\gamma(j,k)\gamma(k,l)\rho(j-l)}\nonumber \\
&\leq c  n^{(2H-1)\vee 0}  \sum_{j,  l=0}^{n-1} \big[(j+1)^{H}- j^{H}\big]\big[(l+1)^{H}- l^{H}\big]\abs{ \rho(j-l)} .\label{2 eistimat of 3 culmulant}
\end{align}
In the above summation, when $j=l$, it is clear that
$$\sum_{j=l} \big[(j+1)^{H}- j^{H}\big]\big[(l+1)^{H}- l^{H}\big]\abs{ \rho(j-l)} = \sum_{j=0}^{n-1} \big[(j+1)^{H}- j^{H}\big]^2\le c n^{(2H-1)\vee 0} .$$
When $j=0<l$, we have 
\begin{align*}\sum_{j=0<l} \big[(j+1)^{H}- j^{H}\big]\big[(l+1)^{H}- l^{H}\big]\abs{ \rho(j-l)} &= \sum_{l=1}^{n-1} \big[(l+1)^{H}- l^{H}\big] \rho(l)\\
&\le c\times  \sum_{l=1}^{n-1} l^{H-1} l^{2H-2}\\ 
& \le c\times n^{ (3H-2)\vee 0},
\end{align*} where in the last line we use Lemma 6.3 of \cite{ND 12}. 
 The symmetry, the inequality \eqref{gamma beta inquality},  the monotonicity of the power function $f(x)=x^{\nu}$ with $x>0,\,\nu<0$, and the change of variable $k=l-j$ imply that the other terms are less than:
\begin{align*}
&2 \times \sum_{ 0<j<l\le n-1} \big[(j+1)^{H}- j^{H}\big]\big[(l+1)^{H}- l^{H}\big]\abs{ \rho(j-l)} \\
&\le c\times \sum_{ 0<j<l\le n-1}  j^{H-1}l^{H-1} (l-j)^{2H-2}\\
&\le  c\times  \sum_{j,k\in \N, j+k<n } j^{2H-2} k^{2H-2}  \le c\times n^{(4H-2)\vee 0}.
\end{align*} where the last line is from Lemma~\ref{beta function inequality}. 
  Plugging the above three estimates into \eqref{2 eistimat of 3 culmulant}, we obtain the inequality  \eqref{2 eiitmate main}.

Finally, we have
\begin{align}
&\abs{\sum_{j, k, l=0}^{n-1}\gamma({j,k})\rho(k-l)\rho(j-l)}\nonumber \\
& \leq \sum_{j, k, l=0}^{n-1}\abs{\gamma({j,k})\rho(k-l)\rho(j-l)}\nonumber \\
&\le \frac{1}{H^2} \sum_{j, k, l=0}^{n-1} \big[(j+1)^{H}- j^{H}\big]\big[(k+1)^{H}- k^{H}\big]\abs{ \rho(k-l)\rho(j-l)} .\label{30 eiitmate main}
\end{align}
In the similar vein, we have that in the summation the contribution of all the terms such that $j=k$ or $k=l$ or $j=l$ or $j=0$ or $k=0$ or $l=0$ are negligible to compare with $n^{(6H-3)\vee 0}$.  The symmetry implies that other terms are less than:
\begin{align*}
&\quad  2 \times \sum_{ 0<j<k< n,\,l\neq j,k} \abs{\gamma({j,k})\rho(k-l)\rho(j-l)}\\
&\le c\times \sum_{ 0<j<k< n,\,l\neq j,k}   j^{H-1}k^{H-1}\abs{ k-l}^{2H-2} \abs{ j-l}^{2H-2} . \end{align*}
According to the distinct orders of $j,\,k,\, l$, we do the change of variables $a=j,\,k-j=b,\,l-k=c$ when $0<j<k<l$, or $a=j,\,l-j=b,\,k-l=c$ when $0<j<l<k$, or $a=l,\,j-l=b,\,k-j=c$ when $0<l<j<k$, and then by the monotonicity of the power function again,  we have
\begin{align*}
& \sum_{ 0<j<k< n,\,l\neq j,k}   j^{H-1}k^{H-1}\abs{ k-l}^{2H-2} \abs{ j-l}^{2H-2}\\
&\le  
3\times \sum_{ a,b,c\in \N, a+b+c<n} a^{2H-2}b^{2H-2}c^{2H-2}\\
&\le c \times n^{(6H-3)\vee 0},
\end{align*}where the last line is from Lemma~\ref{beta function inequality}. 
  Taking the above three inequalities together, we obtain the desired \eqref{3 eiitmate main}.
\end{proof}
In the same way, we can show the following proposition.
\begin{prop}\label{4 main bijiao jichu}
Let $\theta(i,j),\,\gamma(i, j),\,\rho(r )$ be given as in \eqref{gammaij dfn}. When $H\in (0,\,1)$, there exists a positive constant $c$ such that for all $n\ge 1$
\begin{align}
\abs{\sum_{i,j, k, l=0}^{n-1}\gamma({i,j})\gamma(i,k)\gamma(k,l) \gamma(j,l)}& \leq c\times n^{(8H-4)\vee 0}, \nonumber\\
\abs{\sum_{i, j, k, l=0}^{n-1}\gamma(i,j)\gamma(i,k)\gamma(k,l)\rho(j-l)  }&\leq c\times n^{(8H-4)\vee 0}, \label{2 eiitmate main2}\\
\abs{\sum_{i,j, k, l=0}^{n-1}\gamma(i,j)\gamma(i,k)\rho(k-l)\rho(j-l)}&\leq c\times n^{(8H-4)\vee 0}, \label{3 eiitmate main2}\\
\abs{\sum_{i,j, k, l=0}^{n-1}\gamma(i,j)\rho(i-k)\rho(k-l)\rho(j-l)}&\leq c\times n^{(8H-4)\vee 0}. \label{4 eiitmate main2}
\end{align}
\end{prop}
\begin{proof}
 Similarly, 
 we have
\begin{align*}
\abs{\sum_{i,j, k, l=0}^{n-1} \gamma({i,j})\gamma(i,k)\gamma(k,l) \gamma(j,l)}& \leq\sum_{i,j, k, l=0}^{n-1} \abs{\gamma({j,k})\gamma(i,k)\gamma(k,l) \gamma(j,l)}\\
&\le \Big( \frac{1}{H^2} \sum_{j=0}^{n-1}\big[(j+1)^{H}- j^{H}\big]^2\Big)^4\\
&\le c\times n^{(8H-4)\vee 0},
\end{align*}
and 
\begin{align*}
&\abs{\sum_{i, j, k, l=0}^{n-1}\gamma(i,j)\gamma(i,k)\gamma(k,l)\rho(j-l)  }\\
& \leq \sum_{i, j, k, l=0}^{n-1}\abs{\gamma(i,j)\gamma(i,k)\gamma(k,l)\rho(j-l)  }\\
&\le c\times n^{(4H-2)\vee 0} \sum_{j, l=0}^{n-1}\big[(j+1)^{H}- j^{H}\big]\big[(l+1)^{H}- l^{H}\big]\abs{ \rho(j-l)} .\\
&\le c\times n^{(8H-4)\vee 0},
\end{align*} where in the last line we have used the proof of  the inequality  \eqref{2 eiitmate main}, please refer to \eqref{2 eistimat of 3 culmulant}.
Next, we have
\begin{align*}
&\abs{\sum_{i,j, k, l=0}^{n-1}\gamma(i,j)\gamma(i,k)\rho(k-l)\rho(j-l)}\\
&\le \sum_{i,j, k, l=0}^{n-1}\abs{\gamma(i,j)\gamma(i,k)\rho(k-l)\rho(j-l)}\\
&\leq c\times n^{(2H-1)\vee 0} \sum_{j,k, l=0}^{n-1}\big[(j+1)^{H}- j^{H}\big]\big[(l+1)^{H}- l^{H}\big]\abs{\rho(k-l) \rho(j-l)}\\
&\leq c\times n^{(8H-4)\vee 0}
\end{align*} where in the last line we have used the proof of the inequality \eqref{3 eiitmate main}, please refer to \eqref{30 eiitmate main}.

Finally, we have 
\begin{align*}
&\abs{\sum_{i,j, k, l=0}^{n-1}\gamma(i,j)\rho(i-k)\rho(k-l)\rho(j-l)}\\
&\le \sum_{i,j, k, l=0}^{n-1} \abs{\gamma(i,j)\rho(i-k)\rho(k-l)\rho(j-l)}\\
&\leq  \frac{1}{H^2}\sum_{i,j,k, l=0}^{n-1}\big[(i+1)^{H}- i^{H}\big]\big[(j+1)^{H}- j^{H}\big]\abs{\rho(i-k)\rho(k-l) \rho(j-l)}.
\end{align*}
It is easy to show that when any two index of $i,j,k,l$ are equal or any index vanishes, the contribution to the sum are negligible to compare with $n^{(8H-4)\vee 0}$.  The symmetry implies that other terms is less than 
\begin{align*}
4\times \sum_{0<i<j<n,0<k< l<n,\,i,j\neq k,l}  i^{H-1} j^{H-1}\abs{i-k}^{2H-2}\abs{k-l}^{2H-2} \abs{j-l}^{2H-2}.
\end{align*}
According to the distinct orders of $i,\,j,\,k,\, l$, we do the change of variables $a=i,\,j-i=b,\,k-j=c, l-k=d$ when $0<i<j<k<l$, or $a=i,\,k-i=b,\,j-k=c, l-j=d$ when $0<i<k<j<l$, or  $a=j,\,k-i=b,\,l-k=c,\,j-l=d$ when $0<i<k<j<l$, or $a=k,\,i-k=b,\,j-i=c,\,l-j=d$ when $0<k<i<j<l$, or or $a=k,\,i-k=b,\,l-i=c,\,j-l=d$ when $0<k<i<j<l$, or $a=k,\,l-k=b,\,i-l=c,\,j-i=d$ when $0<k<i<j<l$, and then by the monotonicity of the power function again,  we have
\begin{align*}
& \sum_{0<i<j<n,0<k< l<n,\,i,j\neq k,l}  j^{H-1}k^{H-1}\abs{ k-l}^{2H-2} \abs{ j-l}^{2H-2}\\
&\le 6 \sum_{ a,b,c,d\in \N, a+b+c+d<n} a^{2H-2}b^{2H-2}c^{2H-2}d^{2H-2} \\
&\le c \times n^{(8H-4)\vee 0},
\end{align*}where the last line is from Lemma~\ref{beta function inequality}.
\end{proof}
\noindent{\it Proof of Theorem~\ref{main thm 1}.\,} 
We will discuss exclusively the case $H \in(0,\frac34)$ since the case $H= \frac34$ is similar. 
Recall  \eqref{zn exp} and \eqref{gn exp},  the expressions of $Z_n$ and $g_n$. Denote 
\begin{equation*}
F_n:= \frac{Z_n}{\sqrt{n}}=\frac{ I_{2}(g_{n})}{\sqrt{n}}.
\end{equation*} 
First, the identities (6.2-6.3)  of \cite{BBNP 12} imply that 
  \begin{equation}
\kappa_3(F_n)=\E[{F_n}^3]=\frac{8}{n^{\frac32}}{\sum_{j,k,l=0}^{n-1}\theta(j,k)\theta(k,l)\theta(j,l)}
\end{equation}
\begin{equation}\label{gn contraction}
\kappa_4(F_n)={\E[{F_n^4}]-3\E[F_n^2]^2}=\frac{48}{n^2} \norm{g_n\otimes_1 g_n }^2_{\FH^{\otimes 2}}=\frac{48}{n^2}{\sum_{j,k,l=0}^{n-1}\theta(i,j)\theta(j,k)\theta(k,l)\theta(j,l)}
\end{equation}
\\
The symmetry and \eqref{gammaij dfn} imply that 
 \begin{align*}
  {\sum_{j,k,l=0}^{n-1}\theta(j,k)\theta(k,l)\theta(j,l)}&=\sum_{j,k,l=0}^{n-1}\rho(j-k)\rho(k-l)\rho(j-l) +\sum_{j, k, l=0}^{n-1}\gamma(j,k)\gamma(k,l)\gamma(j,l) \\
  &+ 3\sum_{j, k, l=0}^{n-1}\gamma(j,k)\gamma(k,l)\rho(j-l) + 3 \sum_{j, k, l=0}^{n-1}\gamma({j,k})\rho(k-l)\rho(j-l)
 \end{align*}
Rearranging, and using Proposition~\ref{rop 3 bijiao} we have that there exists a positive constant $c$ such that 
\begin{align*}
\abs{\kappa_3(F_n)- \frac{8}{n^{\frac32}} {\sum_{j,k,l=0}^{n-1}\rho(j,k)\rho(k,l)\rho(j,l)}}&\le c\times n^{(6H-3)\vee 0 -\frac32},
\end{align*} which together with Proposotion~\ref{3 jichu bijiao}, implies that when $H\in (0,\,\frac23) $, 
\begin{align}\label{zuihou bijiao 03}
\kappa_3(F_n)\asymp n^{-\frac12};
\end{align} and when $H=\frac23$,\begin{align}
\kappa_3(F_n)\asymp n^{-\frac12}\log ^2 n;
\end{align} 
and when $H\in (\frac23,\frac34) $, 
\begin{align}
\abs{\kappa_3(F_n)}\le n^{\frac12(4H-3)};
\end{align} 

In the same vein, Proposition~\ref{4 main bijiao jichu} implies that there exists a positive constant $c$ such that 
\begin{align*}
\abs{\kappa_4(F_n)- \frac{48}{n^2} \sum_{i,j,k,l=0}^{n-1}\rho(i-j)\rho(i-k)\rho(k-l)\rho(j-l)}&\le c\times n^{(8H-6)\vee 0},
\end{align*}which together with Proposotion~\ref{3 jichu bijiao}, implies that when $H\in (0,\,\frac58) $, 
\begin{align}
\kappa_4(F_n)\asymp n^{-1};
\end{align} and when $H=\frac58$,\begin{align}
\kappa_4(F_n)\asymp n^{-1}\log ^3 n;
\end{align} 
and when $H\in (\frac58,\frac34) $, 
\begin{align}\label{zuihou bijiao}
\abs{\kappa_4(F_n)}\le n^{ 8H-6};
\end{align} 
 Combing \eqref{zuihou bijiao 03}-\eqref{zuihou bijiao} with Theorem~\ref{fm.theorem} and Propositions~\ref{3 jichu bijiao}, \,\ref{sigma2 limit prop}, we obtain the desired result since 
 \begin{align*}
\kappa_4(V_n)=\frac{n^{\frac32}}{\sigma_n^3} \kappa_3(F_n),\qquad \kappa_4(V_n)=\frac{n^2}{\sigma_n^4} \kappa_4(F_n).
 \end{align*}
{\hfill\large{$\Box$}}

\noindent{\it Proof of Theorem~\ref{main thm 2}.\,} 
 From Theorem~\ref{main thm 1}, $V_n$ satisfies the CLT.  Hence,  we need only to check the conditions (1) and (2) of Theorem~\ref{asclt chufad}
 are valid. We will discuss exclusively the case $H \in(0,\frac34)$ since the case $H= \frac34$ is similar. 
  
Recall $V_n=I_2(f_n)$ where
\begin{equation}
f_n=\frac{1}{\sigma_n} g_n=\frac{\sqrt{n}}{\sigma_n} \frac{g_n}{\sqrt{n}},
\end{equation}
which together with \eqref{gn contraction} and Proposition \ref{sigma2 limit prop} implies that
\begin{align*}
\norm{f_n \otimes_1 f_n }^2_{\FH^{\otimes 2}}&\le c\times \frac{1}{n^2} \norm{g_n\otimes_1 g_n }^2_{\FH^{\otimes 2}}\\
&\le c\times  \left\{
      \begin{array}{ll}
 \frac{1}{n}, & \quad \text{if }  H\in (0,\,\frac58),\\
 {\frac{(logn)^3}{n},} & \quad \text{if } H = \frac58, \\
n^{8H-6}, &\quad \text{if } H\in (\frac58,\, \frac34).
  \end{array}
\right.
\end{align*} Hence, the condition (1) of Theorem~\ref{asclt chufad} is valid.

To check the condition (2) of Theorem~\ref{asclt chufad}, first noting that $\E[V_n^2]=2\norm{f_n}^2_{\FH^{\otimes 2}}=1$, we need only to show that when $0<k<l$, 
the following inequality holds:
\begin{align}\label{condition 2 inequality}
\abs{\innp{f_k,f_l}_{\FH^{\otimes 2}}}\le c\times \Big[\sqrt{\frac{k}{l}} +(kl)^{(2H-1)\vee 0 -\frac12}\Big].
 \end{align}
 In fact, we have 
\begin{align}
\abs{\innp{f_k,f_l}_{\FH^{\otimes 2}}}& \le c\times \frac{1}{\sqrt{kl}} \abs{\innp{g_k,g_l}_{\FH^{\otimes 2}}}\nonumber\\
&=c\times \frac{1}{\sqrt{kl}}{\sum_{i=0}^{k-1} \sum_{j=0}^{l-1} \theta^2(i,j) }\nonumber\\
&\le c\times \frac{2}{\sqrt{kl}}{\sum_{i=0}^{k-1} \sum_{j=0}^{l-1} \big[\gamma^2(i,j) +\rho^2(i-j)\big]}. \label{fenjie fkl}
 \end{align}
 \cite[Theorem 5.1]{BNT 10} implies that 
 \begin{align*}
\frac{1}{\sqrt{kl}}{\sum_{i=0}^{k-1} \sum_{j=0}^{l-1} \rho^2(i-j)} \le c\times \sqrt{\frac{k}{l}}.
 \end{align*}
 It follows from the inequality \eqref{gamma estimate} that 
  \begin{align*}
\frac{1}{\sqrt{kl}}{\sum_{i=0}^{k-1} \sum_{j=0}^{l-1} \gamma^2(i,j) }& \le c\times \frac{1}{\sqrt{kl}} \sum_{i=0}^{k-1} [(i+1)^{H}-i^{H}]^2 \times \sum_{j=0}^{l-1}[(j+1)^{H}-j^{H}]^2\\
&\le c\times  (kl)^{(2H-1)\vee 0- \frac12},
 \end{align*}where in the last line we have used the inequality \eqref{gamma estimate 2}.  Plugging the above two inequalities into \eqref{fenjie fkl}, we obtain the desired \eqref{condition 2 inequality}. Hence, the ASCLT holds when  $H\in(0,\frac34)$.
{\hfill\large{$\Box$}}

\section*{Acknowledgements}
Many thanks to Prof. X. Huang for valuable comments on Lemma~\ref{lemm integra by part}. The work of Yong Chen is supported by National Natural Science Foundation of China (No.11961033)

 \end{document}